\documentclass[twoside,10pt,a4paper]{newFNLstyle}
\usepackage{amsmath,graphicx}
\usepackage{cite}

\begin{document}

\volnumpagesyear{0}{0}{000--000}{2009}
\dates{received date}{revised date}{accepted date}

\title{Equilibrium Solution to the Lowest Unique Positive Integer Game}

\authorsone{SEUNG KI BAEK and SEBASTIAN BERNHARDSSON}
\affiliationone{Department of Theoretical Physics, Ume{\aa} University, 901
87 Ume{\aa}, Sweden}
\mailingone{garuda@tp.umu.se, sebbeb@tp.umu.se}
 
\maketitle

\markboth{Equilibrium Solution to the LUPI game}{Baek and Bernhardsson}

\pagestyle{myheadings}


\begin{abstract}
We address the equilibrium concept of a reverse auction game so that
no one can enhance the individual payoff by a unilateral change when all
the others follow a certain strategy. In this approach
the combinatorial possibilities to consider become very much involved even for
a small number of players, which has hindered a precise analysis in previous
works.
We here present a systematic way to reach the solution for a general number of
players, and show that this game is an example of conflict between the group and
the individual interests.
\end{abstract}

\section{Introduction}

Game theory deals with a situation where each player's payoff is dependent
not only on her own behavior but also on other players'. Players will
generally have conflicting interests with each other,
but they may get better off by interacting with others in various
ways~\cite{gintis}.
For example, in the minority game~\cite{minor}, a player gets a point when
her choice is minor among players. Thus, even if a choice happened to be
successful at previous rounds, it is hard to remain minor with that choice
as more and more players will also select it. Although it is not possible
for all the players to win together, the average probability of winning
can become enhanced if they behave in proper ways.

Recently, an extreme version of the minority game was proposed as a special
case of the reverse auction, which is called the {\em lowest unique positive
integer} (LUPI) game ~\cite{zeng}. This game has $n > 2$ players, and each of
them may choose one integer from $1$ to $n$. A player wins a point by
choosing the lowest unique number. That means, simply choosing the number
$1$ is not a good strategy since it is very likely to be chosen by other
players, too.
In Ref.~\cite{flitney}, this LUPI game was analyzed in terms of the
Nash equilibrium (NE): Let us imagine that every player participates
in a game with her own strategy, which is publicly known to the other
players. If no player can improve her payoff by changing her strategy alone,
the set of the players' strategies is called a NE. Even though this is an
adequate framework to analyze a game, the actual enumeration of
combinatorial possibilities becomes hard to manage as the number of players
increases. Although Ref.~\cite{flitney} tried to find an expression for general $n$,
many of the possible cases were missing in the probability calculations, even for $n=4$.
According to the suggested formula in Ref.\ \cite{flitney}, all the loosing bids should be 
above the winning number or all should be below the winning number. Furthermore, if they 
are below the winning number, they should all be on the same number. However, there exist 
many more cases that would generate the same winning number.
We henceforth present a way to take all the possibilities into account correctly for general 
$n$ and analyze the results thereby obtained.

\section{Equilibrium}
\label{sec:ne}

Suppose that each player draws numbers from her own probability
distribution, which we will denote as her strategy.
As explained above, we will find a NE in this strategy
space.
Among several different NE's in the LUPI game,
we will be concerned with one that can be prevalent
among these $n$ players, as studied in Ref.~\cite{flitney}. This approach
actually corresponds to the stability concept
in the evolutionary game theory~\cite{weibull}.
Let us imagine that $n-1$ people are using a certain
strategy, $\vec{p} = (p_1,p_2, \ldots p_n)$, where $p_i$ means the
probability to choose number $i$.
The idea is to choose $\vec{p}$ so that the chance of winning for the $n$th
player becomes independent of her strategy.
We consider all the possible cases that those $n-1$
people can make. These can be symbolically represented by
\begin{equation}
Z_0 = \left(\sum_{i=1}^n p_i \right)^{n-1},
\label{eq:total}
\end{equation}
from which all the combinations can be obtained as separate terms.
For example, in a three-players game we get
$Z_0=(p_1+p_2+p_3)^2 = p_1^2+p_2^2+p_3^3+2p_1p_2+2p_1p_3+2p_2p_3$,
where the first three terms on the right-hand side mean that players 1 and 2 chose the same number 1, 2 or 3, respectively.
The other three terms correspond to the cases when players 1 or 2 have chosen two different numbers.
Together, these terms represent all the possible outcomes of the game (as seen from the third player) since they add up to the probability 1 due to the normalization condition $\sum_i^n p_i =1$.

Now, starting from Eq.\ (\ref{eq:total}), we extract all the cases where there is a winner, among the $n-1$ players, at number 1.
These are characterized by all the terms which contain only one $p_1$.
In other words, if we write Eq.~(\ref{eq:total}) as a polynomial in $p_1$:
\begin{equation*}
Z_0 = A(p_2,\ldots,p_n) +
B(p_2,\ldots,p_n) p_1 + C(p_2,\ldots,p_n) p_1^2 + \ldots,
\end{equation*}
those cases are expressed by the second term and
one can easily find it by $B = \left. d Z_0 / d p_1 \right|_{p_1=0}$.
In other words, we find the cases with no winner at number 1 as
\begin{eqnarray*}
Z_1 &=& Z_0 - \left.\frac{d Z_0}{d p_1}\right|_{p_1=0} p_1\nonumber\\
&=& \left(\sum_{i=1}^n p_i \right)^{n-1} - (n-1)p_1 \left(\sum_{i \neq 1}^n
p_i \right)^{n-2} .
\end{eqnarray*}
If we are further to exclude the cases where there is a winner at number
2, we apply the same operation on $Z_1$ as follows:
\begin{eqnarray*}
Z_2 &=& Z_1 - \left.\frac{d Z_1}{d p_2}\right|_{p_2=0} p_2\nonumber\\
&=& \left(\sum_{i=1}^n p_i \right)^{n-1} - (n-1)p_1 \left(\sum_{i \neq 1}^n
p_i \right)^{n-2}\nonumber\\
&-& (n-1)p_2 \left(\sum_{i \neq 2}^n p_i \right)^{n-2}
+ (n-1)(n-2) p_1 p_2 \left(\sum_{i \neq 1,2}^n p_i \right)^{n-3}.
\end{eqnarray*}
Generalizing this idea, we can find cases where there is no
winner up to an arbitrary number $i$ by the following recursion relation:
\begin{equation}
Z_i = Z_{i-1} - \left.\frac{d Z_{i-1}}{d p_i}\right|_{p_i=0} p_i.
\label{eq:recur}
\end{equation}
Formally, we may define $L_i$ as the {\em linearity projection operator}
with respect to $p_i$. Let us denote a generic polynomial of
$p_1,\ldots,p_n$ as $Q$, which may have an index to indicate more than one
polynomial.
The operator $L_i$ can be algebraically represented as $L_i [Q] =
p_i E_i \cdot D_i [Q]$, where $E_i$ is elimination of $p_i$ by substituting
zero and $D_i$ is differentiation with respect to $p_i$.
In addition, $a$ and $b$ mean some coefficients independent of $p_i$. We
then have the following relations for these operators:
\begin{itemize}
\item $L_i[a p_i^n] = a p_i \delta_{1,n}$
\item $L_i \cdot L_i [Q] = L_i [Q]$
\item $L_i \cdot L_j [Q] = L_j \cdot L_i [Q]$
\item $L_i [a Q_1 + b Q_2] = a L_i[Q_1] + b L_i[Q_2]$
\item $E_j \cdot L_i [Q] = L_i \cdot E_j [Q]$
\item $L_i [Z_j] = 0$, ~~~~if $i \le j$
\end{itemize}
where $\delta_{1,n}$ is the Kronecker delta.
It immediately follows that
\begin{equation}
Z_k = \left\{ \prod_{i=1}^k \left(1-L_i \right) \right\} [Z_0].
\label{eq:op}
\end{equation}
Let us now calculate the probability for the $n$th player to win the LUPI game by
choosing number $i$. Then within the other $n-1$ people, there should be
no winner up to $i-1$ and no one should choose the number $i$.
The last condition can be easily imposed by the same trick. That is, we
simply substitute zero for $p_i$ to exclude every case where $p_i$ appears.
Therefore, the probability of our concern is simply written as $c_i
(\vec{p}) = E_i [Z_{i-1}] = \left\{ \prod_{j=1}^{i-1} \left(1-L_j \right)
\right\} \cdot E_i [Z_0] $ and the $n$th player's expected payoff is
\begin{equation*}
W(\vec{\pi};\vec{p}) = \sum_{i=1}^{n} c_i(\vec{p}) \pi_i,
\end{equation*}
where $\vec{\pi}=(\pi_1,\pi_2,...,\pi_n)$ is the $n$th player's strategy.

\begin{figure}
\centering{
\includegraphics[width=0.49\textwidth]{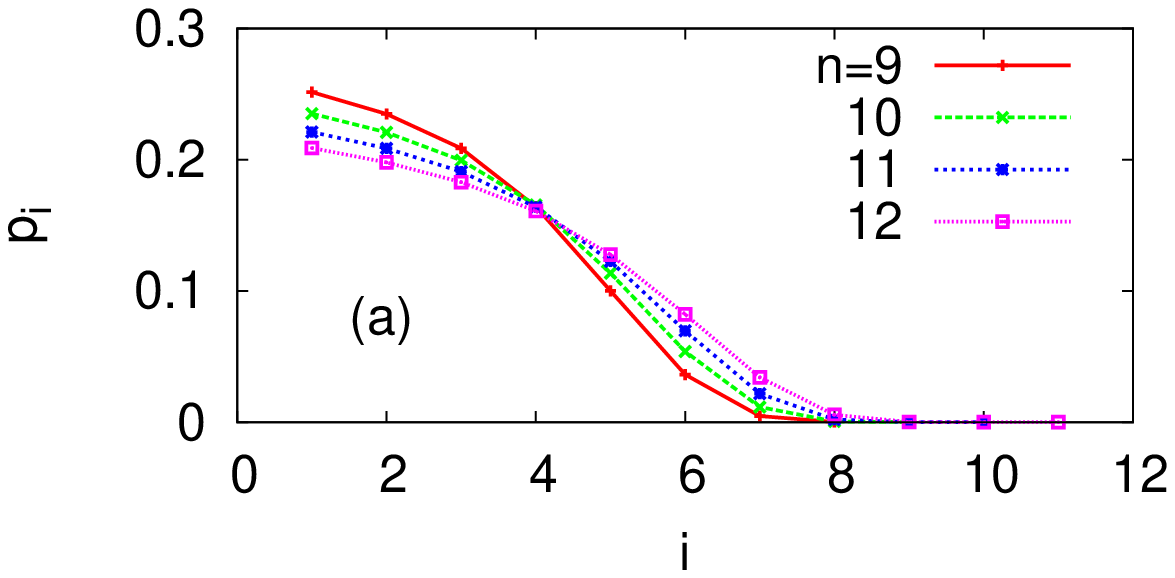}
\includegraphics[width=0.49\textwidth]{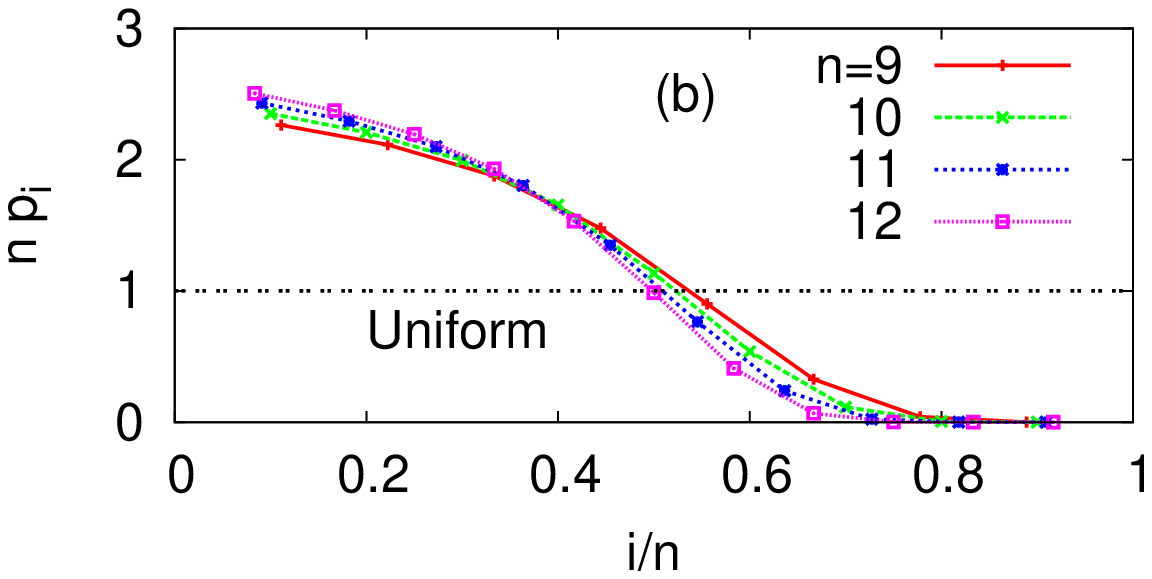}
}
\caption{(a) NE strategies, $\vec{p}_{\rm NE}$, with varying the number of
players, $n$. (b) Scaled plots with respect to $n$. As $n$ increase, the
solution deviates more from the uniform solution (dotted).}
\label{fig:ne}
\end{figure}

For example,
\begin{eqnarray}
c_1 &=& \left. Z_0 \right|_{p_1=0} = \left( \sum_{i \neq 1}^n p_i
\right)^{n-1} = (1-p_1)^{n-1},\label{eq:c1}\\
c_2 &=& \left. Z_1 \right|_{p_2=0} = \left(\sum_{i \neq 2}^n p_i
\right)^{n-1} - (n-1)p_1 \left(\sum_{i \neq 1,2}^n p_i
\right)^{n-2}\nonumber\\
&=& (1-p_2)^{n-1} - (n-1)p_1 (1-p_1-p_2)^{n-2},\label{eq:c2}\\
c_3 &=& \left. Z_2 \right|_{p_3=0}\nonumber\\
&=& \left(\sum_{i \neq 3}^n p_i \right)^{n-1}
- (n-1)p_1 \left(\sum_{i \neq 1,3}^n p_i \right)^{n-2}\nonumber\\
&-& (n-1)p_1 \left(\sum_{i \neq 2,3}^n p_i \right)^{n-2}
+ (n-1)(n-2)p_1 p_2 \left(\sum_{i \neq 1,2,3}^n p_i \right)^{n-3}\nonumber\\
&=& (1-p_3)^{n-1} - (n-1)p_1 (1-p_1-p_3)^{n-2}\nonumber\\
&-& (n-1)p_2 (1-p_2-p_3)^{n-2}\nonumber\\
&+& (n-1)(n-2) p_1 p_2 (1-p_1-p_2-p_3)^{n-3}.\label{eq:c3}
\end{eqnarray}
From the normalization condition, $p_n = 1 - \sum_{i=1}^{n-1} p_i$, indeed
we have only $n-1$ degrees of freedom in choosing $\vec{p}$ so that
\begin{equation*}
W(\vec{\pi};\vec{p}) = \sum_{i=1}^{n-1} (c_i - c_n) \pi_i + c_n.
\end{equation*}
The NE solution, $\vec{p}_{\rm NE}$, should satisfy $c_i
(\vec{p}_{\rm NE}) = $ constant for every $i$ so that $W$ cannot be better by
changing $\vec{\pi}$. That is, the chance of winning is the same on all numbers. For $n=3$, all the above calculation coincides with
that presented in Ref.~\cite{flitney}.
Note that we have $n-1$ degrees of freedom and $n-1$ equations so we can
solve these simultaneously. Figure~\ref{fig:ne}(a) shows solutions for some
$n$ values, obtained by using the Newton method. In
Fig.~\ref{fig:ne}(b), we plot its scaled version, considering that the
horizontal axis naturally scales with $n$, and the vertical axis roughly
with $1/n$.

\begin{figure}
\centering{
\includegraphics[width=0.49\textwidth]{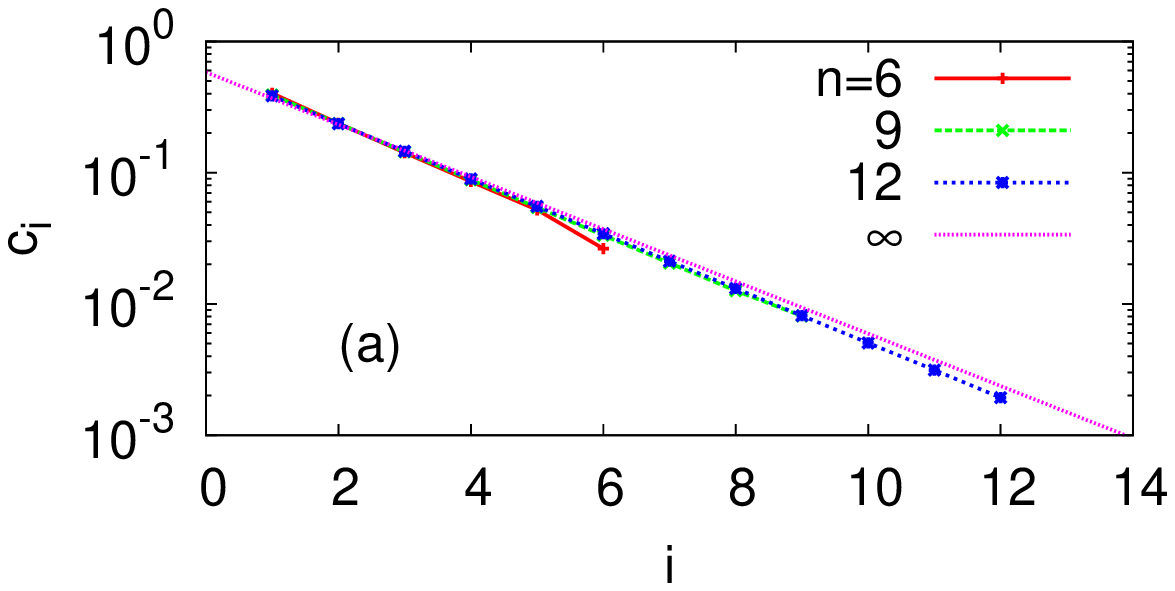}
\includegraphics[width=0.49\textwidth]{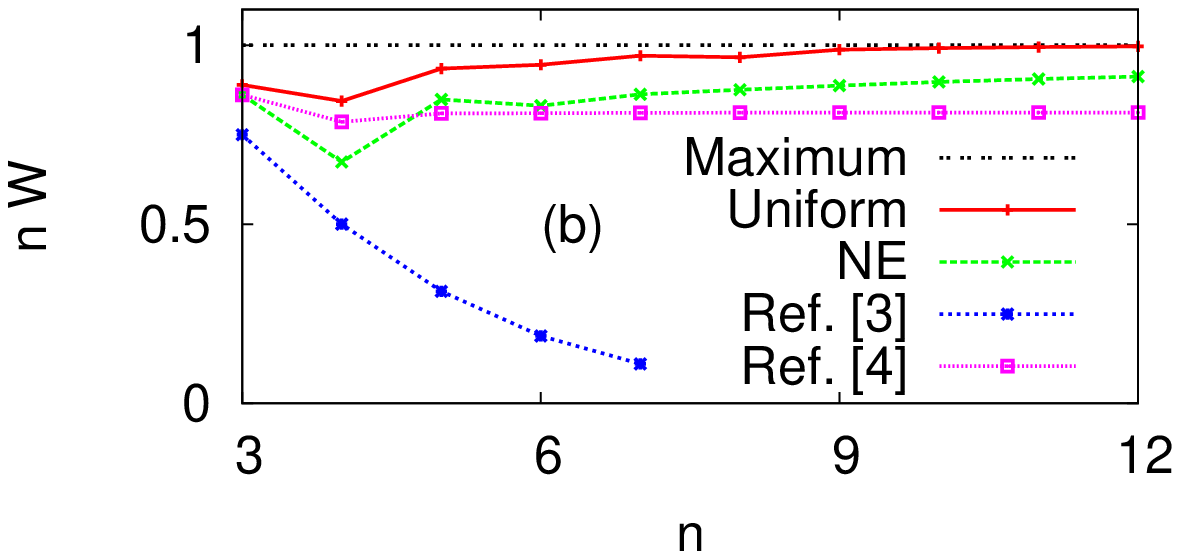}}
\caption{(a) The expected chance of winning by choosing $i$, when the other
$n-1$ players follow the uniform solution defined as $p_i = 1/n$.
(b) Comparison of outcomes scaled with $n$. The uniform solution
gives an outcome that rapidly converges to the theoretical maximum
$W=1/n$ (dotted), while the NE outcome remains below that. The other two
curves represent outcomes from the strategies suggested in previous works.}
\label{fig:bd}
\end{figure}

Setting $c_1 = c_2$ from Eqs.~(\ref{eq:c1}) and (\ref{eq:c2}),
we find
\begin{equation}
(1-p_2)^{n-1} - (1-p_1)^{n-1} = (n-1) p_1 (1-p_1-p_2)^{n-2}.
\label{eq:unequal}
\end{equation}
Since the right-hand side should be positive, we can conclude that $p_2 <
p_1$ for every finite $n$, which means the equilibrium solution cannot be
uniform.
Its implication is remarkable: Let us imagine that everyone employs the same
strategy $\vec{p}$. In other words, $\vec{\pi} = \vec{p}$.
The question is what the strategy $\vec{p}$ should be
in order to maximize the expected payoff, $W(\vec{p};\vec{p})$.
For example, for $n=3$, we find
\begin{equation*}
W(\vec{p};\vec{p}) = (1-p_1) (1-p_2) (1-p_3),
\end{equation*}
which is maximized at a uniform solution, $\vec{p}_{\ast} =
(1/3,1/3,1/3)$. This is clearly different from the NE solution,
$\vec{p}_{\rm NE} = (2\sqrt{3}-3, 2-\sqrt{3}, 2-\sqrt{3})$~\cite{flitney}.
The expression is not so simple for $n>3$, but solving $\nabla_{\vec{p}}
W(\vec{p};\vec{p}) = 0$ under $\sum_i p_i = 1$, we still see that the
uniform solution gives the best outcome to the population on average.
Although we confirmed this only for $n \le 12$,
this result is plausible since the risk of overlapping
choices will be minimized by the uniform solution.
On the other hand, we already know that the equilibrium strategy cannot be
uniform. Indeed, if all the people follow $\vec{p}_{\ast} = (1/n,
\ldots, 1/n)$, a player's chance of winning rapidly decays with her
chosen number $i$ [Fig.~\ref{fig:bd}(a)]. For example,
one may use the formulation above to obtain
$c_i = e^{-1}(1-e^{-1})^{i-1}$
for large $n$, which coincide with the multinomial analysis in
Ref.~\cite{zeng}.
Therefore, each player is motivated to
deviate from this uniform solution by choosing low numbers more
frequently.
In short, there exists a dilemma between the {\em best} and the
{\em equilibrium} strategy.
Note that
since each play only has one winner, the theoretical upper limit for the
gain per play is $1/n$. Let us then compare these three cases:  the
equilibrium outcome, the best possible outcome, and the theoretical upper
limit [Fig~\ref{fig:bd}(b)]. While the uniform solution quickly converges to
the theoretical maximum, the NE solution remains suboptimal in this
plot. For comparison, we also plot outcomes from the strategies in
previous works: Ref.~\cite{zeng} suggested a strategy with $p_1=p_2=1/2$
which gives $W(\vec{p}, \vec{p}) = 2^{1-n}$, while in Ref.~\cite{flitney}
the solution is approximately written as $p_i = 2^{-i}$ for $i<n$ with
$p_n = 2^{1-n}$. Note that the latter one actually performs better than
presented in Ref.~\cite{flitney}, since we have taken all the possibilities
into consideration. In any case, these become smaller
than the outcome from our NE solution.

As Eq.~(\ref{eq:op}) indicates, the number of terms to
consider increases exponentially as $n$ grows, making the exact enumeration
of probabilities intractable. Although our current computational
resources allow us to arrive only up to $n=12$, our procedure provides a
systematic way to take all the combinatorial cases into account for any $n$, in principle.
We believe that the behaviors shown in Fig.~\ref{fig:bd} give a clue to
anticipate $\vec{p}_{\rm NE}$ in the limit of large $n$. Specifically, we
sketch a way to perform this task: The NE solution is obtained when $c_i=c_j$ for all $i$ and $j$.
This means that we can set $c_i$ to a certain constant $c_0$ for $i=1,\ldots,n$. The point is that $c_i$ has only one more
variable, $p_i$, than $c_{i-1}$, so that we can attack each variable $p_i$
one by one. For example, Eq.~(\ref{eq:c1}) with $c_1 = c_0$ gives $p_1 (c_0) = 1 -
c_0^{1/(n-1)}$. Substituting this value into Eq.~(\ref{eq:c2}) and with $c_2=c_0$ gives a value for $p_2(c_0)$, and so on.
If the correct NE-value is chosen so that $c_0=c_{\rm NE}$ then the equilibrium solution, $\vec{p}_{\rm NE}$ is obtained.
When there exist several real solutions for $p_i$ then the correct root can be found using the restriction $0 < p_i < 1$ and the normalization condition $\sum_{j=1}^i p_j < 1$.
The procedure for a three-players ($n=3$) game, where we have $c_0 = c_{\rm NE} = 28-16\sqrt{3} \approx 0.287$ (Fig.~\ref{fig:bd}b), goes as follows:

\begin{eqnarray}
 &c_1& = c_0 = (1-p_1)^2 \Rightarrow p_1 \approx 0.4643.\nonumber\\
 &c_2& = c_0 = (1-p_2)^2 - 2p_1(1-p_1-p_2) \nonumber\\
&\Rightarrow& [\textrm{inserting } p_1] \Rightarrow p_2 \approx 0.2684.\nonumber\\
&c_3& = c_0 = (1-p_3)^2 - 2p_1(1-p_1-p_3) - 2p_2(1-p_2-p_3) + 2p_1p_2\nonumber\\
&\Rightarrow& [\textrm{inserting } p_1 \textrm{ and } p_2] \Rightarrow p_3 \approx 0.2673 \pm 0.0246i.\nonumber
\end{eqnarray}

The fact that the last solution is complex indicates that the value of $c_0$ is not exactly correct. Adding more digits decreases the imaginary part while the real part moves closer to the correct solution (e.g.,\ $c_0 = 0.287187$ gives $p_3 = 0.267949 - 5.04978\times 10^{-4}i$).
Since there is no general formula for roots of polynomials of degree larger than four, it seems that this procedure should be carried out numerically. Figure 3 shows the behavior of $c_i$ for $i=1$ to $4$ and $n=9$. The solution of $p_i$ is obtained in each step by finding the crossing point $c_i=c_0$ (where $p_{i-1}$ has been inserted into $c_i$).

\begin{figure}
\centering{
\includegraphics[width=0.8\textwidth]{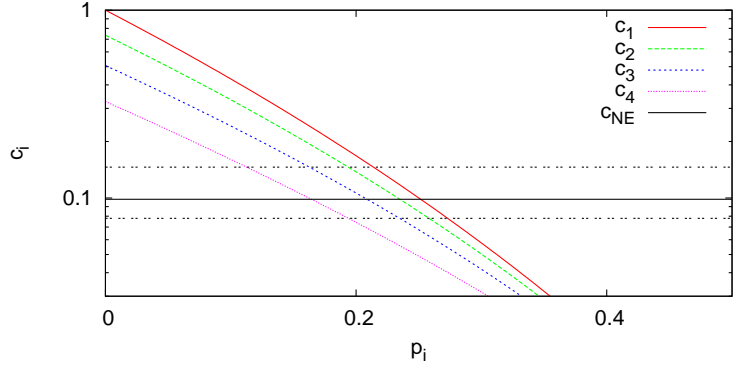}
\caption{The expressions $c_i$ as functions of $p_i$ for $n=9$ players. A Nash equilibrium solution is obtained when $c_i=c_{\rm NE}$ for all $i \leq n$. $p_1$ is found by the condition $c_1=c_{\rm NE}$. Inserting this value into $c_2$ allows for a determination of  $p_2$ as the crossing point between $c_2$ and $c_{\rm NE}$, and so on. Or the value of $c_{\rm NE}$ can be obtained in this way by finding $c_0$ which fulfills the normalization condition $\sum_{i=1}^n p_i=1$. The two dotted horizontal lines represent the interval within which $c_{\rm NE}$ must lie if we use the expressions up to $c_4$ instead of $c_n$.}}
\label{fig:3}
\end{figure}

To sum up, the benefit of this procedure is that if we know the correct NE-value, $c_0=c_{\rm NE}$, then we can calculate the probabilities one by one up to any $i$ we want or are capable of doing within our computational resources, for any $n$. In many situations it might, for example, be enough to know the solution up to certain $i$ above which $p_i$'s are negligible. The method of solving for all $i$ simultaneously instead results in an \emph{all-or-nothing} situation.
Of course, the down side is that we need to know $c_{\rm NE}$ in advance. 
However, recalling that $c_{\rm NE} = W_{\rm NE} \equiv W (\vec{p}_{\rm NE}; \vec{p}_{\rm NE})$, we see the
limiting behavior $\lim_{n \rightarrow \infty} n c_{\rm NE} = \lim_{n \rightarrow
\infty} n W_{\rm NE}$ in Fig.~\ref{fig:bd}(b). Putting such a $c_0=c_{\rm NE}$ into this
procedure may thus yield $\vec{p}_{\rm NE}$ for large $n$.
The solution obtained in this way up to $i=4$ for $n=9$ and $c_0 = 0.0985 \approx c_{\rm NE}$ (see Fig.\ 3) is $\vec{p}=(0.2515,0.2349,0.2087,0.1643,...)$. This can be compared to the result obtained from solving for all $i$ simultaneously, which is $\vec{p}=(0.2515,0.2348,0.2086,0.1641,...)$.

This method could also be used to find the actual NE solution, without knowing $c_{\rm NE}$, by iterating the procedure until a self-consistent solution is found for all $i$. If this is again beyond our computational power, one could at least determine an interval in which $c_0=c_{\rm NE}$ lies. Let's say that we use expressions for $c_i$ up to $i=j$. If $c_0$ is too small then the normalization gets violated by $\sum_{i=1}^j p_i > 1$. If $c_0$ is too large then the crossing points result in probabilities that are too small to add up to one in total, as long as $p_i$ is a monotonically decreasing function of $i$. Thus, the largest possible sum of all $n$ probabilities is obtained if $p_i$ is uniform for $i\ge j$, so the sum $\sum_{i=1}^j p_i + (n-j)p_j$ must be larger or equal to one.
The interval obtained for $n=9$ and $j=4$ is marked by the two horizontal, dotted lines in Fig.\ 3 ($c_{0, min}=0.078$ and $c_{0, max}=0.146$). This interval becomes narrower as higher values of $i$ are addressed. Note that when $c_0$ is changed all the curves $c_i$ for $i > 1$ also change since we get new $p_j$'s with $j < i$ to be inserted into $c_i$.

\section{Discussion}
\label{sec:discuss}

We have shown a way of taking all the combinatorial possibilities into 
account for the LUPI (lowest unique positive integer) game and
proposed a procedure to find a Nash equilibrium for a general number of players, $n$.
It turns out, however, that the number of terms to consider grows exponentially in $n$,
and that the computational capability of our personal computer can only solve the problem 
up to $n=12$ players. To deal with this problem we also suggest an alternative way of 
solving for the equilibrium strategy  
for each value of $i$ at a time, by inserting an estimate of the average expected payoff of each player.
This procedure allows the problem to be solved sequentially instead of solving the whole problem simultaneously.

We also have found that a uniform solution distributes nearly optimal
outcomes equally among all the players. However, this would be driven to
a suboptimal equilibrium solution, so we may regard this $n$-person game as
an example of a social dilemma to manifest a conflict between group and
individual interests.

\section*{Acknowledgments}
S.K.B. acknowledges the support from the Swedish Research Council
with the Grant No. 621-2002-4135.


\end{document}